\documentclass[11pt,reqno]{amsart}\usepackage[letterpaper, portrait, margin=1in, headheight=0pt]{geometry}
\usepackage{mathtools}
\usepackage{enumerate}
\usepackage{amsmath,amsthm,amssymb,amsfonts}
\usepackage{stmaryrd}
\usepackage{hyperref}
\usepackage[normalem]{ulem}
\usepackage{soul}
\usepackage{xcolor}
\usepackage{tikz}
\usepackage{tikz-cd}
\usetikzlibrary{calc, positioning, fit, shapes.misc,arrows,bending}
\usepackage{subcaption}
\usepackage{blkarray}
\usepackage{cleveref}
\newtheorem{thm}{Theorem} 
\theoremstyle{definition}
\newtheorem{exam}[thm]{Example}
\newtheorem{defi}[thm]{Definition}

\newtheorem{cor}[thm]{Corollary}

\newtheorem{prop}[thm]{Proposition}

\makeatletter
\DeclareRobustCommand{\cev}[1]{%
  \mathpalette\do@cev{#1}%
}
\newcommand{\do@cev}[2]{%
  \fix@cev{#1}{+}%
  \reflectbox{$\m@th#1\vec{\reflectbox{$\fix@cev{#1}{-}\m@th#1#2\fix@cev{#1}{+}$}}$}%
  \fix@cev{#1}{-}%
}
\newcommand{\fix@cev}[2]{%
  \ifx#1\displaystyle
    \mkern#23mu
  \else
    \ifx#1\textstyle
      \mkern#23mu
    \else
      \ifx#1\scriptstyle
        \mkern#22mu
      \else
        \mkern#22mu
      \fi
    \fi
  \fi
}

\makeatother

\newcommand{\dvecEG}{%
  \mathrel{\vbox{\offinterlineskip\ialign{%
    ##\hfil\cr
    $\scriptscriptstyle\hspace{.2ex}\leftrightarrows$\cr
    \noalign{\kern-.1ex}
    $E(G)$\cr
    }}}}

\usepackage{titlesec}

\titleformat{\section}
       {\normalfont
       \fontsize{14}{17}\bfseries}{\thesection}{1em}{}

\titleformat{\subsection}
       {\normalfont\fontsize{12}{17}\bfseries}{\thesubsection}{1em}{}

\title{Some orientation theorems for restricted DP-colorings of graphs}

\begin{document}

\setstcolor{red}
\maketitle
\begin{center} Ian Gossett\

\end{center}

\begin{abstract}

We define $Z$\textit{-signable} correspondence assignments on multigraphs, which generalize \textit{good} correspondence assignments as introduced by Kaul and Mudrock. 
We introduce an auxiliary digraph that allows us to prove an Alon-Tarsi style theorem for DP-colorings from $Z$-signable correspondence assignments on multigraphs, and apply this theorem to obtain three DP-coloring analogs of the Alon-Tarsi theorem for arbitrary correspondence assignments as corollaries. We illustrate the use of these corollaries for DP-colorings on a restricted class of correspondence assignments on toroidal grids.
        
\end{abstract}

\section{Introduction}
\vspace{.2cm}

DP-colorings of graphs (originally called correspondence colorings) are a generalization of list colorings that were introduced by Dvorak and Postle in \cite{MR3758240} and used as a tool in their proof that planar graphs without cycles of length $4$ to $8$ are $3$-choosable. DP-colorings have since become a popular topic of study, with one of the overarching research goals being to characterize the similarities and differences between DP-colorings and list colorings. Though some results for list colorings hold for DP-colorings in general, many results are known not to carry over to DP-colorings (see, e.g., \cite{MR4262023},\cite{MR3518419},\cite{MR3889157},\cite{MR3758240}). In particular, it is well-known that the Alon-Tarsi list coloring theorem does not apply `as is' to DP-colorings, in the sense that we cannot replace list coloring with DP-coloring in the statement of the theorem.

In this work, we define \textit{generalized signable} and $Z$\textit{-signable} correspondence assignments of multigraphs. These classes generalize \textit{good} correspondence assignments, originally introduced by Kaul and Mudrock in \cite{MR4137634}, where they were called \textit{good prime covers}. DP-colorings from generalized signable and $Z$-signable correspondence assignments specialize to colorings of signed graphs as defined in \cite{MR3484719} and \cite{MR0675866}, $\mathbb{Z}_p$-colorings of signed graphs as defined in \cite{MR4003656}, and list colorings of signed graphs as defined in \cite{MR3802151}. Furthermore, DP-colorings from these classes are closely related to certain subclasses of field colorings \cite{MR4366879}, $A$-colorings \cite{MR1186753}, and colorings of gain graphs \cite{MR1328292}.

 Our main result is \Cref{thm:Eulerianaux}, in which we use an auxiliary digraph construction to prove an Alon-Tarsi style theorem for DP-colorings from $Z$-signable correspondence assignments on multigraphs. 
 As applications of the main result, we obtain three DP-coloring analogs of the Alon-Tarsi theorem (Corollaries \ref{cor:multizsign}, \ref{cor:multisign},\ref{cor:multigood}). This is achieved by viewing an arbitrary correspondence assignment on 
 a simple graph as a Z-signable correspondence assignment on an appropriate multigraph, instead. Since no other theorems exist in the literature which use an auxiliary digraph construction to obtain orientation theorems for DP-colorings, these theorems are the first of their kind. 

The paper proceeds as follows. In Section 2, we establish basic notational conventions and give necessary background on list colorings, orientations, and DP-colorings. In Section 3.1, we introduce the main classes of correspondence assignments that will be studied, and in Section 3.2, we highlight the relationship between these classes and the Combinatorial Nullstellensatz. In Section 4.1, we introduce an auxiliary digraph construction that will be used in our main result,  in Section 4.2, we prove the main result, and in Section 4.3, we prove some immediate corollaries for closely related classes of of correspondence assignments. In Section 5, we present three DP-coloring analogs of the Alon-Tarsi list coloring theorem for simple graphs and give an application to certain restricted DP-colorings of toroidal grids.  In Section 6, we give some concluding remarks.

  \vspace{.2cm}

\section{Preliminaries}\label{sec:background}

\subsection{Background and notation}
In this section, we provide prerequisite background concerning list colorings, orientations, and DP-colorings. In general, by a \textit{graph}  we will mean a finite, loopless multigraph and by \textit{digraph} (or \textit{directed graph}) we mean a finite, loopless multidigraph, though we will often be explicit about this for the sake of clarity. By \textit{loopless}, we mean that no vertex is self adjacent. If $G=(V,E)$ is a graph, and $e\in E$, we write $e= \{v,w\}$ to specify that $v$ and $w$ are the endpoints of $e$.  (This is a slight abuse of notation, since the endpoints do not determine $e$ uniquely.) Similarly, when $D=(V,\vec{E})$ is a digraph, we write $\vec{e}=(v,w)$ to express that $\vec{e}$ is a directed edge from $v$ to $w$. If $G=(V,E)$ is a graph such that $e_1=\{u,v\}\in E$ and $e_2=\{u,v\}\in E$ implies that $e_1=e_2$, then $G$ is called simple.  Similarly, a digraph $D=(V,\vec{E})$ is simple if $\vec{e}_1=(v,w)\in \vec{E}$ and $\vec{e}_2=(v,w)\in \vec{E}$ implies $\vec{e}_1=\vec{e}_2$.

Let $D=(V,\vec{E})$ be a digraph. For each $v\in V(D)$, let $d^+_D(v)$ denote the out-degree of $v$ in $D$, and let $d^-_D(v)$ denote the in-degree of $v$ in $D$. An \textit{orientation} of a simple graph $G=(V,E)$ is a directed graph $D=(V,\vec{E})$ such that $V(D)=V(G)$, and $\{u,v\}\in E(G)$ if and only if exactly one of $(u,v)$ or $(v,u)$ is in  $\vec{E}$. Orientations of multigraphs are defined analogously. If $D$ is an orientation of $G$ and $e$ is an edge of $G$, we will denote the directed edge in $D$ corresponding to $e$ by $\vec{e}$. 
The following prerequisite definitions and theorems that will be used in this paper. 

\begin{defi}
    Let $G$ be graph and $L=\{L(v)\}_{v\in V(G)}$ be an assignment of lists to the vertices of $G$. If there exists a function $f:V\rightarrow \bigcup_{v\in V}{L(v)} $ such that $f(v)\in L(v)$ for each $v\in V$, and $f(u)\neq f(v)$ for all $\{u,v\}\in E(G)$, then $f$ is called an $\pmb{L}$\textbf{-coloring} of $G$. The $\textbf{list chromatic number}$ of $G$, denoted $\chi_{\ell}(G)$, is the least integer $k$ such that for any assignment of lists $L=\{L(v)\}_{v\in V}$ with $|L(v)|\geq k$ for each $v\in V(G)$, $G$ is $L$-colorable. 
\end{defi}

\begin{defi}
Let $D$ be a digraph. $D$ is said to be \textbf{Eulerian} if for each $v\in V(D)$, $d^+_D(v)=d^-_D(v)$. (Note that Eulerian digraphs are not required to be connected.)
\end{defi}

\begin{defi}
A digraph is called $\textbf{even}$ if it has an even number of edges, and $\textbf{odd}$ if it has an odd number of edges. 
Let $EE(D)$ denote the number of even spanning Eulerian subdigraphs of $D$, and $EO(D)$ denote the number of odd spanning Eulerian subdigraphs of $D$. 
\end{defi}

The following theorem is known as \textit{Combinatorial Nullstellensatz}. It is the main algebraic tool used in this work.

\begin{thm}\label{thm:CN}(Alon, \cite{MR1684621})
Let $F$ be a field and $f$ a polynomial in $F[x_1,x_2,...,x_n]$. Let $d_1,...,d_n$ be nonnegative integers such that the total degree of $f$ is $\sum_{i=1}^{n}d_i$. For each $1\leq i \leq n$, let $L_i$ be a set of at least $d
_i+1$ elements of $F$. If the coefficient in $f$ of $\prod_{i=1}^{n}x^{d_i}_i$ is nonzero, then there exists $x\in L_1\times L_2\times \cdots \times L_n$ such that $f(x)\neq 0$.
\end{thm}

The following theorem is an application of  \Cref{thm:CN} that is known as the Alon-Tarsi theorem. 

\begin{thm} \label{thm:at}(Alon-Tarsi, \cite{MR1179249}) Let $G=(V,E)$ be a graph and let $D$ be an orientation of $G$. If  $EE(D)\neq EO(D)$, and if $L(v)$ is a set of $d^+_D(v)+1$ distinct integers for each $v\in V(G)$, then there is a proper vertex coloring $c:V\rightarrow \mathbb{Z}$ such that $c(v)\in L(v)$ for each $v\in V$. 
\end{thm}

\subsection{DP-colorings}

We will now turn our attention to DP-colorings and related terminology. The following definition of a correspondence coloring (since renamed DP-coloring) was first presented by Dvorak and Postle in \cite{MR3758240} for simple graphs. The definition was later extended to multigraphs in \cite{MR3686937}. 

\begin{defi} \label{def:cor}
  Let $G$ be a multigraph. \vspace{.2cm} 

\begin{itemize}

\item A \textbf{correspondence assignment} for $G$ consists of an assignment of lists $L=\{L(v)\}_{v\in V}$ and a function $C$ that to every edge $e=\{u,v\}$ assigns a partial matching $C_e$ between $\{u\}\times L(u)$ and $\{v\}\times L(v)$.\vspace{.2cm} 

\item An $(L,C)$-coloring of $G$ is a function $f$ that to each vertex $v\in V(G)$ assigns a color $f(v)\in L(v)$, such that for every $e=\{u,v\}\in E(G)$, the vertices $(u,f(u))$ and $(v,f(v))$ are non-adjacent in $C_e$. If such an $(L,C)$ coloring exists, we say that $G$ is $(L,C)$-colorable. 

\end{itemize}

 Colorings of the above type are called \textbf{DP-colorings}. 
\end{defi}

%

The following diagram shows an example of a correspondence assignment $(L,C)$ on $C_4$, with $L(v)=\{1,2\}$ for each $v\in V(C_4)$. 
 As has become the custom, we represent each list $L(v)$ ``inside" its corresponding vertex, and draw the matchings for each edge. (The reader will quickly ascertain that $C_4$ is not $(L,C)$-colorable.) 

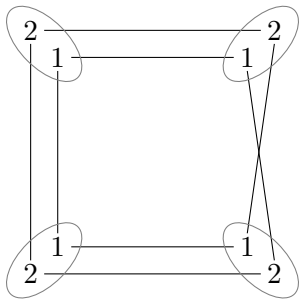
\begin{figure}[h]

\vspace{.1cm}
\label{fig:badcor} \begin{center}
\begin{tikzpicture}[scale=1.8]
\node at (1.75,.5){} ;



\node at (-.2,-.2){1};
\node at (-.4,-.4){2};
\node at (-.2,1.2){1};
\node at (-.4,1.4){2};
\node at (1.2,1.2){1};
\node at (1.4,1.4){2};
\node at (1.2,-.2){1};
\node at (1.4,-.4){2};

\draw[] (-.2,-.1)--(-.2,1.1);
\draw[] (-.4,-.3)--(-.4,1.3);
\draw[] (-.1,1.2)--(1.1,1.2);
\draw[] (-.3,1.4)--(1.3,1.4);
\draw[] (-.1,-.2)--(1.1,-.2);
\draw[] (-.3,-.4)--(1.3,-.4);

\draw[] (1.2,1.1)--(1.4,-.3);
\draw[] (1.4,1.3)--(1.2,-.1);


the\draw[rotate around={45:(-.3,1.3)},gray] (-.3,1.3) ellipse (5pt and 10pt);
\draw[rotate around={-45:(1.3,1.3)},gray] (1.3,1.3) ellipse (5pt and 10pt);
\draw[rotate around={-45:(-.3,-.3)},gray] (-.3,-.3) ellipse (5pt and 10pt);
\draw[rotate around={45:(1.3,-.3)},gray] (1.3,-.3) ellipse (5pt and 10pt);
\end{tikzpicture}
\end{center}
\end{figure}

List colorings can be seen as the special case of DP-colorings where each partial matching $C_{\{u,v\}}$ is given by $\{(u,c_1),(v,c_2)\}\in C_{\{u,v\}}$ if and only if $c_1=c_2$. Following the terminology in \cite{MR3758240}, if $(L,C)$ is a correspondence on $G$ and $e=\{u,v\}\in E(G)$ is such that $\{(u,c_1),(v,c_2)\}\in C_{\{u,v\}}$ implies $c_1=c_2$, then $e$ is said to be \textit{straight}.

 \Cref{def:cor} gives rise to a natural notion of the \textit{DP-chromatic number of a multigraph \cite{MR3686937}}: 

\begin{defi} \label{def:chidp}
The \textbf{DP-chromatic number} (also called the correspondence chromatic number) of a multigraph $G$, denoted $\chi_{DP}(G)$, is the least $k$ such that for every correspondence assignment $(L,C)$ on $G$ with $|L(v)|\geq k$ for each $v\in V$, there exists an $(L,C)$-coloring of $G$. 
\end{defi}

In \cite{MR3758240}, a notion of equivalent correspondence assignments was introduced.
 \begin{defi} Two correspondence assignments $(L,C)$ and $(L',C')$ of a graph $G$ are \textbf{equivalent} if there exists a family of bijections $\{h_v\}_{v\in V}$, with $h_v:L(v)\rightarrow L'(v)$ for each $v$, such that for each $e=\{v,w\}\in E(G)$, $\{(v,c_1),(w,c_2)\}\in C_e$ if and only if $\{(v,h(c_1),(w,h(c_2)\}\in C'_e$.
\end{defi}

The bijections $h_v$ in the above definition can be seen as a `renaming' of the list elements of $(L,C)$, while preserving the structure of the matchings. It is immediate that if a graph $G$ is $(L,C)$-colorable, then $G$ is $(L',C')$-colorable for all equivalent correspondence assignments $(L',C')$.  

Our goal in this work is to relate orientations and DP-colorability. We start by examining the simple example of $C_4$.  The correspondence assignment shown in \Cref{fig:badcor} demonstrates that the Alon-Tarsi theorem does not apply if we replace list colorings with DP-colorings in its statement; a cyclic orientation $D$ on $C_4$ has $EE(D)\neq EO(D)$, but the correspondence assignment shown in \Cref{fig:badcor} proves that even if each vertex $v\in V(C_4)$ is equipped with a list of $d^+_D(v)+1=2$ elements, this is not enough to guarantee a DP-coloring from an arbitrary correspondence assignment on these lists. 
 
However, if we instead consider the digraph $D$ shown below, 
whose underlying simple undirected graph is $C_4$, then it is indeed true that any correspondence $(L,C)$  on $C_4$ with $|L(v)|\geq d^+_{D}(v)+1$, admits an $(L,C)$-coloring of $C_4$. 
This basic idea of `doubling up' edges, when necessary, will be generalized and made precise in Corollaries \ref{cor:multizsign},\ref{cor:multisign}, and \ref{cor:multigood}. 

\begin{figure}[ht]

\vspace{.1cm}
\begin{center}
\begin{tikzpicture}[scale=2.5]

\tikzset{vertex/.style = {shape=circle,fill=black, draw, inner sep=0pt,minimum size=2.5mm}}
\tikzset{edge/.style = {-{Stealth[scale=1.8]}
}}
\node[vertex] (a) at  (0,0) {};
\node[vertex] (b) at  (0,1) {};
\node[vertex] (c) at  (1,1) {};
\node[vertex] (d) at  (1,0) {};
\draw[edge] (a) to (b);
\draw[edge] (b) to (c);
\draw[edge] (c) to[bend left] (d);
\draw[edge] (d) to[bend left] (c);
\draw[edge] (d) to (a);

\end{tikzpicture}
\end{center} 
\end{figure}

\section{Restricted DP-colorings and the Combinatorial Nullstellensatz}
In this section, we define \textit{good} (originally defined in \cite{MR4137634}), \textit{signable}, \textit{generalized signable}, and $Z$\textit{-signable} correspondence assignments, and highlight their relation to the Combinatorial Nullstellensatz.

\subsection{Classes of F-correspondence assignments} 

We will consider correspondence assignments whose list elements come from some field, $F$. We remark, however, that the list elements in any correspondence assignment can be renamed to come from some field. (See the discussion in Section 6.) 

\begin{defi} Let $F$ be a field. A correspondence assignment $(L,C)$ on a graph $G$ is called an $\pmb{F}$\textbf{-correspondence assignment} if $L(v)\subseteq F$ for each $v\in V$.  
\end{defi}

 The following notion of a \textit{good} $F$-correspondence assignment was introduced in \cite{MR4137634} (though, formulated in a different way), where the study of the relation between DP-colorings and the Combinatorial Nullstellensatz was initiated. The original definition was only for simple graphs, but the same definition works for $F$-correspondence assignments on multigraphs:
 
\begin{defi}\label{def:good}
Let $G=(V,E)$ be a multigraph, and fix an orientation $D=(V,\vec{E})$ of $G$. Let $(L,C)$ be an $F$-correspondence assignment on $G$. Suppose that $e=\{v,w\}\in E(G)$ is such that $\vec{e}=(v,w)\in \vec{E}(D)$. We say that $e$ is \textbf{good} with respect to $(L,C)$ if there exists some $a_{\vec{e}}\in F$  such that  $c_1-c_2=a_{\vec{e}}$ whenever $\{(v,c_1),(w,c_2)\}\in C_e$. $(L,C)$ is said to be \textbf{good} if every $e\in E(G)$ is good. 
\end{defi}

Note that if $e$ is a straight edge, we can always set $a_{\vec{e}}=0$. Note also that the definition of a good edge is independent of the choice of the orientation $D$, since if $\vec{e}=(v,w)$ and $\cev{e}=(w,v)$, we can set $a_{\cev{e}}=-a_{\vec{e}}$.

The following definition of \textit{signable} correspondence assignments generalizes good correspondence assignments.

\begin{defi}\label{def:signable}
    Let $G=(V,E)$ be a multigraph, and fix an orientation $D$ of $G$. Let $(L,C)$ be an $F$-correspondence assignment on $G$. If there exists a sign function $\sigma_C:E(G)\rightarrow \{1,-1\}$ so that for each $e=\{v,w\}\in E(G)$, with $\vec{e}=(v,w)$, there exists some $a_{\vec{e}}\in F$ such that  $c_1-\sigma_C(e)c_2=a_{\vec{e}}$ whenever and $\{(v,c_1),(w,c_2)\}\in C_{\{v,w\}}$, then we call $(L,C)$ a \textbf{signable} correspondence assignment on $G$. 
\end{defi}

Hence, if $e\in E(G)$ and $(L,C)$ is a signable correspondence assignment on $G$ such that $\sigma_C(e)=1$, then $e$ is good, and if $\sigma_C(e)=-1$ we will call $e$ \textbf{bad} (following the terminology from \cite{MR4137634}). A good correspondence assignment can therefore be viewed as a signable correspondence assignment such that $\sigma_C(e)=1$ for all $e\in E(G)$. 

We briefly discuss the relationship between signable correspondence assignments and signed colorings of signed graphs. A \textit{signed graph} is a pair $(G,\sigma)$, where $G$ is a graph and $\sigma:E(G)\rightarrow \{-1,1\}$.  A (signed) coloring in $k$-colors of $(G,\sigma)$, as defined in \cite{MR0675866}, is a function $\psi:V(G)\rightarrow \{-k,-k+1,...-1,0,1,...,k-1,k\}$ such that for each $e=\{v,w\}\in E(G)$, we have that $\psi(v)\neq \sigma(e)\psi(w)$. In \cite{JIN2016234}, the idea of signed colorings was extended to signed list colorings. Given a signed graph, $(G,\sigma)$, and an assignment of lists $L=\{L(v)\}_{v\in V}$, with $L(v)\subseteq \mathbb{Z}$ for each $v$, an $(L,\sigma)$-coloring is a function $\psi:V(G)\rightarrow \bigcup_{v\in V}L(v)$ such that $\psi(v)\in L(v)$ for each $v$, and for each $e=\{v,w\}\in E(G)$, we have that $\psi(v)\neq \sigma(e)\psi(w)$. If $(G,\sigma)$ admits an $(L,\sigma)$-coloring, we say $(G,\sigma)$ is $L$-colorable. Furthermore, signed $\mathbb{Z}_p$ colorings, where $L(v)$ is instead a subset of $ \mathbb{Z}_p$ for each $v$, were considered in \cite{MR4003656}.


If we begin with a signed graph, $(G,\sigma)$, then  DP- colorings from the set of signable $F$-correspondence assignments $(L,C)$ on $G$ which also admit $\sigma$ as a sign function are a generalization of signed colorings (\cite{MR3484719},\cite{MR0675866}), 
signed $\mathbb{Z}_p$-colorings \cite{MR4003656}, and signed list colorings \cite{MR3802151} of the signed graph $(G,\sigma)$. However, due to the allowance of partial matchings in  \Cref{def:cor}, there exist correspondence assignments with sign function $\sigma$ which admit DP-colorings that are not also signed colorings of $(G,\sigma)$. For example, consider the following $\mathbb{R}$-correspondence assignment $(L,C)$ on $K_2$. 

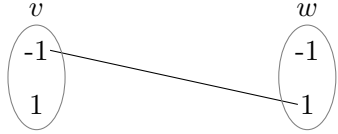
\begin{center}
\begin{tikzpicture}[scale=1.8]
\node at (1.75,.5){} ;

\node at (-1,.5){$v$};
\node at (-1,.2){-1};
\node at (-1,-.2){1};

\node at (1,.5){$w$};
\node at (1,.2){-1};
\node at (1,-.2){1};

\draw[] (-.9,.2)--(.93,-.2);

\draw[gray] (-1,0) ellipse (6pt and 11pt);
\draw[gray] (1,0) ellipse (6pt and 11pt);

\end{tikzpicture}
\end{center}

Evidently, $(L,C)$ admits the sign function $\sigma$ given by $\sigma(\{v,w\})=-1$. However, the function $f$ defined by $f(v)=1$ and $f(w)=-1$ is an $(L,C)$ coloring that is not also a signed coloring of the signed graph $(K_2,\sigma)$, because $f(v)=\sigma(\{v,w\})f(w)$.

We now define \textit{generalized signable} correspondence assignments, a further generalization of signable correspondence assignments. 

\begin{defi}\label{def:gensignable} Let $G=(V,E)$ be a multigraph, and fix an orientation $D=(V,\vec{E})$ of $G$. Suppose that $F$ is a field and $(L,C)$ is an $F$-correspondence assignment on $G$. A \textbf{generalized sign function} for $(L,C)$ is a function 
 $\varphi_C:E(G)\rightarrow F\setminus \{0\} $ with the property that for each $\vec{e}=(v,w)\in \vec{E}(D)$, there exists some $a_{\vec{e}}\in F$ such that $c_1-\varphi_C(e)c_2=a_{\vec{e}}$ whenever $\{(v,c_1),(w,c_2)\}\in C_{\{v,w\}}$. If such a function $\varphi_C$ exists, we say that $(L,C)$ is a \textbf{generalized signable} $F$-correspondence assignment on $G$. 
\end{defi}

Thus, a signable correspondence assignment is a generalized signable correspondence assignment such that $\varphi_C(e)\in\{1,-1\}$ for each $e\in E(G)$. The definition of a generalized signable correspondence assignment does not depend on the chosen orientation $D$ of $G$: Suppose that $D'$ is an orientation of $G$ and $\cev{e}=(w,v)\in \vec{E}(D')$, with $\vec{e}=(v,w)\in E(D)$. Then since $c_1-\varphi_C(e)c_2=a_{\vec{e}}$ implies that $c_2-\varphi(e)^{-1}c_1=-\varphi(e)^{-1}a_{\vec{e}}$, we see that $D'$ admits a generalized sign function $\hat{\varphi}_C$ such that $\hat{\varphi}_C(e)=\varphi_C(e)$ for each edge $e$ of $G$ that is oriented in the same direction in $D$ and $D'$, and $\hat{\varphi_C}(e)=\varphi_C(e)^{-1}$ (and $a_{\cev{e}}=-\varphi(e)^{-1}a_{\vec{e}}$) for each edge $e$ of $G$ whose orientation in $D$ is opposite its orientation in $D'$. Hence, the definition of a generalized signable correspondence assignment on $G$ does not depend on the given orientation $D$ of $G$, but the particular map $\varphi_C$ does depend on $D$. When we need to reference this fact explicitly, we will say that $\varphi_C$ is the generalized sign function of $(L,C)$ with respect to $D$.

 The following restriction of generalized signable correspondence assignments will be of particular interest in Sections 4 and 5. 

\begin{defi}\label{def:zsign}
Let $G$ be a graph, and fix an orientation $D$ of $G$. Let $F$ be a field, and let $(L,C)$ be an $F$-correspondence of $G$. Let $\langle 1 \rangle$ denote the additive subgroup of $F$ that is generated by the element $1$. If there exists a generalized sign function  $\varphi_C:E(G)\rightarrow \langle 1 \rangle \setminus \{0\}$ such that $c_1-\varphi_C(e)c_2=a_{\vec{e}}$ whenever $\{(v,c_1),(w,c_2)\}\in C_{\{v,w\}}$ and $(v,w)\in \vec{E}(D)$, then we say that $(L,C)$ is $\pmb{Z}$\textbf{-signable} with respect to $D$.
\end{defi}

Note that, since $\varphi_C(e)^{-1}$ may not be in $\langle 1 \rangle\setminus \{0\}$, the definition of $Z$-signable is \textit{not} independent of the orientation $D$, so we must always specify that $(L,C)$ is $Z$-signable with respect to the orientation $D$. Note also that if $F=\mathbb{Z}_p=\{0,1,\ldots, p-1\}$, with $p$ prime (and with arithmetic mod $p$), then the definitions of generalized signable and $Z$-signable coincide. 

For Z-signable $F$-correspondence assignments where $F$ is a field of characteristic $0$ or $F=\mathbb{Z}_p$, we will find it useful to consider both a sign function $\sigma:E(G)\rightarrow \{1,-1\}$ and a 
``positive" function $\varphi^+_C:E(G)\rightarrow \langle 1 \rangle\setminus \{0\}$, such that  for each edge $e$ of $G$ we have $\varphi_C(e)=\sigma_C(e)\varphi^+_C(e)$, and $\varphi^+_C(e)>0$. In the case that $F=\mathbb{Z}_p$, $\sigma$ can be any sign function as long as we pick $\varphi^+_C(e)$ accordingly. In the the case that $F$ is a field of characteristic $0$, we must take $\sigma_C(e)=-1$ when $\varphi_C(e)<0$, and set $\varphi^+_C(e)=-\varphi_C(e)$, and we must set $\sigma_C(e)=1$ when $\varphi_C(e)>0$ and set $\varphi^+_C(e)=\varphi_C(e)$. (Here, we are identifying $\langle 1 \rangle$ with $\mathbb{Z}$ in the natural way.)

\subsection{Applications of the Combinatorial Nullstellensatz}

The following proposition is a straightforward generalization of results from both \cite{MR4137634} and \cite{MR4003656} to generalized signable correspondence assignments on multigraphs. We note also that the idea for the proof of \Cref{prop:genkm} is contained in \cite{MR4366879}. We include a proof here for the sake of completeness, since the proposition has not yet been stated in the context of generalized signable correspondence assignments on multigraphs.

\begin{prop}\label{prop:genkm}
Let $G=(V,E)$ be a multigraph with $V=\{1,2,...,n\}$, and let $D$ be an orientation of $G$. Let $F$ be a field and $(L,C)$ a generalized signable $F$-correspondence assignment on $G$, and let $\varphi_C:E(G)\rightarrow F \setminus \{0\}$ be the generalized sign function for $(L,C)$ with respect to $D$. Define $h_D\in F[x_1,x_2,...,x_n]$ by $$h_D(x_1,x_2,...,x_n)= \prod_{\vec{e}=(v,w)\in \vec{E}(D)}(x_v-\varphi_C(e)(x_w)).$$

If there exists a monomial $M=\prod_{i=1}^nx^{t_i}_i$ with nonzero coefficient in $h_D$ such that $|L(i)|\geq t_i+1$ for each $i\in V$, then $G$ is $(L,C)$-colorable.
\end{prop}
 \begin{proof}   For each $\vec{e}\in \vec{E}(D)$, define $a_{\vec{e}}$ as in \Cref{def:gensignable}, and define the polynomial $\hat{h}_D\in F[x_1,x_2,...,x_n] $ by  $$\hat{h}_D(x_1,x_2,...x_n)=\prod_{\vec{e}=(v,w)\in \vec{E}(D)}(x_v-\varphi_C(e)x_w-a_{\vec{e}}).$$ Note that if $x\in \prod_{i=1}^{n}L(i)$ is such that $\hat{h}_D(x)\neq 0_F$, then $x$ corresponds to an $(L,C)$-coloring of $G$. Note also that if $M=\prod_{i=1}^nx^{t_i}_i$ such that $\sum_{i=1}^{n}t_i=|E(D)|$, then $M$ has maximum total degree in $\hat{h}$. Furthermore, the coefficient of $M$ in $\hat{h}_D$ is the same as the coefficient of $M$ in $h_D$. Thus, if the coefficient of $M$ in $h_D$ is nonzero, then the coefficient of $M$ in $\hat{h}_D$ is nonzero, and by the Combinatorial Nullstellensatz there must be some $x\in \prod_{i=1}^{n}{L(i)}$ such that $\hat{h}_D(x)\neq 0$. Hence, there exists an $(L,C)$ coloring of $G$.    
\end{proof}

 In \cite{MR4137634}, Kaul and Mudrock made the observation that every $\mathbb{Z}_3$-correspondence assignment is signable, and they proved the following corollary of \Cref{prop:genkm} for the special case of $\mathbb{Z}_3$-correspondence assignments on simple graphs. In \cite{MR4003656}, a similar result to \Cref{cor:3km} for signed colorings was also proved.

\begin{cor}\label{cor:3km}
Let $G=(V,E)$ be a multigraph with $V=\{1,2,...,n\}$, and let $D$ be an orientation of $G$. Let $F$ be an arbitrary field and $(L,C)$ a signable $F$-correspondence assignment on $G$, with sign function $\sigma_C:E(G)\rightarrow \{1,-1\}$. Define $g_D\in F[x_1,x_2,...,x_n]$ by $$g_D(x_1,x_2,...,x_n)= \prod_{\vec{e}=(v,w)\in \vec{E}(D)}(x_v-\sigma_C(e)(x_w)).$$

If there exists a monomial $M=\prod_{i=1}^nx^{t_i}_i$ with nonzero coefficient in $g_D$ such that $|L(v)|\geq t_i+1$ for each $v\in V$, then $G$ is $(L,C)$-colorable.
\end{cor}

 In \cite{MR4003656},  a combinatorial interpretation of the coefficients of the polynomial in \Cref{cor:3km}, in terms of signed Eulerian subdigraphs, was given. In Section 4, we give a combinatorial interpretation in terms of Eulerian subdigraphs of an auxiliary digraph, instead (\Cref{cor:signaux}).

The following corollary of \Cref{prop:genkm} extends a theorem originally proved in \cite{MR4137634} to multigraphs. 

\begin{cor}\label{cor:nullmulti} Let $G=(V,E)$ be a multigraph with $V=\{1,2,...,n\}$, and let $D$ be an orientation of $G$. Let $F$ be an arbitrary field and suppose that $(L,C)$ is a good $F$-correspondence assignment on $G$.  Define $f_D\in F[x_1,x_2,...,x_n]$  by $$f _D(x_1,x_2,...x_n)=\prod_{\vec{e}=(v,w)\in \vec{E}(D)}(x_v-x_w).$$ If there is some monomial $M=\prod_{i=1}^nx^{t_i}_i$ with nonzero coefficient in $f_D$ such that  $|L(i)|\geq t_i+1$ for each $i$, then $G$ is $(L,C)$-colorable.

\end{cor}

\section{Orientations and Z-signable correspondence assignments}
In this section, we prove an Alon-Tarsi style theorem for DP-colorings from $Z$-signable $F$-correspondence assignments, by way of an auxiliary digraph. In order to establish a meaningful combinatorial interpretation of the desired coefficients in the polynomials related to DP-colorings from $F$-correspondence assignments, we restrict to the cases where  $F$ is a field of characteristic $0$ (and identify $\langle 1\rangle$ with $\mathbb{Z}$ in the natural way) or $F=\mathbb{Z}_p$. 

\subsection{An auxiliary digraph construction}
\begin{defi}\label{def:dphi}
Let $G=(V,E)$ be a multigraph, and $D$ an orientation of $G$. Suppose that $F$ is a field of characteristic $0$ or $F=\mathbb{Z}_p$, and that $(L,C)$ is an $F$-correspondence assignment on $G$ that is $Z$-signable with respect to $D$. Let $\varphi_C:E(G)\rightarrow \langle 1 \rangle \setminus \{0\}$ be the generalized sign function of $(L,C)$ with respect to $D$, and write $\varphi_C(e)=\sigma_C(e)\varphi^+_C(e)$ for each $e\in E(G)$, as in the discussion after \Cref{def:zsign}. 

Define the digraph $D_{\sigma,\varphi}$ as follows:\vspace{.2cm}

If $\vec{e}=(v,w)\in \vec{E}(D)$ and $\sigma_C(e)=1$, replace the edge $\vec{e}=(v,w)$ with the following directed graph: 

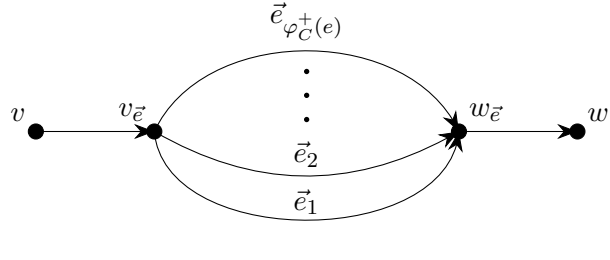
\begin{center}
    \begin{tikzpicture}[scale=.9]

\tikzset{vertex/.style = {shape=circle,fill=black, draw, inner sep=0pt,minimum size=2 mm}}
\tikzset{edge/.style = {-{Stealth[scale=1.6]}
}}


\node[vertex]  at  (-2.5,0) {};
\node[vertex]  at  (-.75,0) {};
\node[vertex]  at (3.75,0){};
\node[vertex]  at (5.5,0){};

\node[anchor=south east] at  (-2.5,0) {$v$};
\node[anchor=south east] at  (-.75,0) {$v_{\vec{e}}$};
\node[anchor=north ] at  (1.5,-.7) {${\vec{e}}_1$};
\node[anchor=north] at  (1.5,0) {${\vec{e}}_2$};

\node[anchor=south west] at  (3.75,0) {$w_{\vec{e}}$};
\node[anchor=south] at  (1.5,1.2) {$\vec{e}_{\varphi^+_C(e)}$};
\node[anchor=south west] at  (5.5,0) {$w$};

\draw[edge] (-2.5,0) to (-.75,0);
\draw[edge] (-.75,0) to [bend right =30] (3.75,0);
\draw[edge] (-.75,0) to [bend left=65] (3.75,0);
\draw[edge] (-.75,0) to  [bend right = 85] (3.75,0);

\draw[edge] (3.75,0) to (5.5,0){};

 \path (1.5,-.3) -- (1.5,1.5) node [font=\huge, midway, sloped] {$\dots$};

\end{tikzpicture};
\end{center}

If $\sigma_C(e)=-1$, replace the edge $\vec{e}=(v,w)$ with the following directed graph:
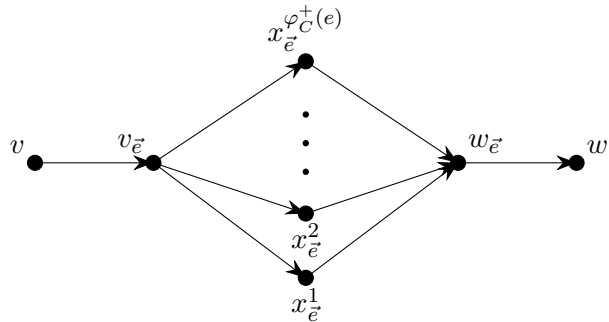
\begin{center}
    \begin{tikzpicture}[scale=.9]

\tikzset{vertex/.style = {shape=circle,fill=black, draw, inner sep=0pt,minimum size=2 mm}}
\tikzset{edge/.style = {-{Stealth[scale=1.6]}
}}


\node[vertex]  at  (-2.5,0) {};
\node[vertex]  at  (-.75,0) {};
\node[vertex]  at  (1.5,1.5) {};
\node[vertex]  at  (1.5,-.75) {};
\node[vertex]  at  (1.5,-1.7) {};
\node[vertex]  at (3.75,0){};
\node[vertex]  at (5.5,0){};

\node[anchor=south east] at  (-2.5,0) {$v$};
\node[anchor=south east] at  (-.75,0) {$v_{\vec{e}}$};
\node[anchor=north ] at  (1.5,-1.7) {$x^1_{\vec{e}}$};
\node[anchor=north] at  (1.5,-.75) {$x^2_{\vec{e}}$};
\node[anchor=south west] at  (3.75,0) {$w_{\vec{e}}$};
\node[anchor=south] at  (1.5,1.5) {$x^{\varphi^+_C(e)}_{\vec{e}}$};
\node[anchor=south west] at  (5.5,0) {$w$};

\draw[edge] (-2.5,0) to (-.75,0);
\draw[edge] (-.75,0) to  (1.5,1.5);
\draw[edge] (-.75,0) to  (1.5,-.75);
\draw[edge] (-.75,0) to  (1.5,-1.7);

\draw[edge] (1.5,1.5) to  (3.75,0);
\draw[edge] (1.5,-.75) to (3.75,0);
\draw[edge] (1.5,-1.7) to  (3.75,0);

\draw[edge] (3.75,0) to (5.5,0){};

 \path (1.5,-.75) -- (1.5,1.5) node [font=\Huge, midway, sloped] {$\dots$};

\end{tikzpicture}
\end{center}

More formally, $D_{\sigma,\varphi}$ is the directed graph with vertex set $$V(D_{\sigma,\varphi})=V(D)\cup \{v_{\vec{e}},w_{\vec{e}}:\vec{e}=(v,w)\in \vec{E}(D)\}\cup \{x^i_{\vec{e}}:\vec{e}\in E(D) \text{ and }\sigma_C(e)=-1, 1\leq i \leq \varphi^+_C(e)\},$$

and edge set \begin{align*}
&\{(v,v_{\vec{e}}),(w_{\vec{e}},w):\vec{e}=(v,w)\in \vec{E}(D)\}\\
\cup&\{\vec{e}_{i}:\sigma_C(e)=1, 1\leq i \leq \varphi^+_C(e)\}\\
\cup &\{(v_{\vec{e}},x^i_{\vec{e}}), (x^i_{\vec{e}},w_{\vec{e}}): \sigma_C(e)=-1, 1\leq i \leq \varphi^+_C(e)\}.\\
\end{align*}
\end{defi}
In the above construction of $D_{\sigma,\varphi}$, each edge $\vec{e}=(v,w)\in \vec{E}(D)$, is replaced by a digraph that has exactly $\varphi^+_C(e)$ distinct directed paths from $v$ to $w$. We call these paths $\gamma$\textit{-paths}:

\begin{defi}\label{def:gamma}
    Let $G$ be a multigraph, $F$ a field of characteristic $0$ or $F=\mathbb{Z}_p$, and let $(L,C)$ be an $F$-correspondence assignment on $G$ that is a $Z$-signable with respect to $D$, with sign function $\varphi_C:E(G)\rightarrow \langle 1 \rangle \setminus \{0\}$. Let $D_{\sigma,\varphi}$ be defined as above. For each $\vec{e}\in E(D)$, and each $1\leq i\leq \varphi^+_C(e)$, define $P_{\vec{e},i}$ by 

\[ P_{\vec{e},i}=\begin{cases} 
(v,v_{\vec{e}}),\vec{e}_{i},(w_{\vec{e}},w) & \text{ if } \sigma_C(e)=1\\
(v,v_{\vec{e}}),(v_{\vec{e}},x^i_{\vec{e}}), (x^i_{\vec{e}},w_{\vec{e}}),(w_{\vec{e}},w) & \text{ if } \sigma_C(e)=-1\\
\end{cases}
\]

The paths $P_{\vec{e},i}$ are called $\pmb{\gamma}$\textbf{-paths}.
\end{defi}

The following proposition will be useful. 

\begin{prop}\label{prop:gamma}  Let $G$ be a multigraph, and let $D$ be an orientation of $G$. Suppose that $F$ is a field of characteristic $0$ or $F=\mathbb{Z}_p$, and let $(L,C)$ be an $F$-correspondence assignment on $G$ that is $Z$-signable with respect to $D$, with sign function $\varphi_C:E(G)\rightarrow \langle 1 \rangle \setminus \{0\}$. A spanning subdigraph $S$ of $D_{\sigma,\varphi}$ is Eulerian if and only if $d^+(v)=d^-(v)$ for each $v\in V(S)\cap V(D)$ and the edge set of $S$ is a (possibly empty) union of pairwise edge-disjoint $\gamma$-paths. 
\end{prop}
\begin{proof} Suppose $S$ is a spanning Eulerian subdigraph of $D_{\sigma,\varphi}$. Then, by definition, we must have that $d^+(v)=d^-(v)$ for each $v\in V(S)\cap V(D)$. If $u\in V(S)\setminus V(D)$, then $u$ is an internal vertex of some $\gamma$-path in $D_{\sigma,\varphi}$, by construction. Furthermore, $d^-_{D_{\sigma,\varphi}}(u)=1$, which implies that we must have either $d^-_{S}(u)=0$ or $d^-_{S}(u)=1$. If $d^-_{S}(u)=0$, then $d^+_S(u)=0$, and no edges incident to $u$ contribute to $\vec{E}(S)$. If $d^-_S(u)=1$, then $d^+_S(u)=1$, and it follows that the edge set of $S$ is composed of disjoint paths  whose internal vertices are all vertices in $ V(S)\setminus V(D)$, and whose end vertices must be in $V(S)\cap V(D)$.  Evidently, the only such paths in $D_{\sigma,\varphi}$ are the $\gamma$-paths, so this proves the forward direction.  

For the reverse implication, since $d^+_S(v)=d^-_S(v)$ for all $v\in V(S)\cap V(D)$, it remains only to prove that that $d^+_S(u)=d^-_S(u)$ for each $u\in V(S)\setminus V(D)$. As noted above, each such $u$ can only be an internal vertex of some $\gamma$-path in  $D_{\sigma,\varphi}$. Furthermore, by the construction of $D_{\sigma,\varphi}$, if two $\gamma$-paths of $D_{\sigma,\varphi}$ are edge-disjoint, then they have no internal vertices in common. Thus, since the $\gamma$-paths in the union are assumed to be edge-disjoint, if $u$ lies on some $\gamma$-path whose edges are included in $\vec{E}(S)$,  it lies on exactly one such $\gamma$-path, and we must have that $d^+_S(u)=d^-_S(u)=1$. If $u$ does not lie on one of the $\gamma$-paths in the union, then $d^+_S(u)=d^-_S(u)=0$. Hence, $d^+_S(v)=d^-_S(v)$ for all $v \in V(S)$, and $S$ is Eulerian. 
\end{proof}
\subsection{An orientation theorem for Z-signable correspondence assignments}
The proof of Theorem \ref{thm:Eulerianaux} is modeled off of a proof of the Alon-Tarsi theorem given in \cite{MR4748249}.
If $F$ is a field of characteristic $0$, we identify $\langle 1 \rangle$ (the additive subgroup generated by $1$) with $\mathbb{Z}$ in the natural way. If $n\in\mathbb{Z}$ and $F=\mathbb{Z}_p$, we write $n\equiv 0_F$ if $n \equiv 0 \text{ (mod }p)$. If $F$ has characteristic $0$, we write $n\equiv 0_F$ if $n=0$. 

\begin{thm}\label{thm:Eulerianaux}Let $G$ be a multigraph, and let $F$ be a field of characteristic $0$ or $F=\mathbb{Z}_p$. Suppose that $(L,C)$ is an $F$-correspondence assignment on $G$ that is $Z$-signable with respect to some orientation $D$, such that $EE(D_{\sigma,\varphi})-EO(D_{\sigma,\varphi})\not\equiv 0_F$ and $|L(v)|\geq d^+_D(v)+1$ for each $v\in V$. Then $G$ is $(L,C)$-colorable. 

\end{thm}

\begin{proof}
Let us first consider the case where $F$ has characteristic $0$. Write $V(G)=\{1,2,...,n\}$, and let $h_D\in F[x_1,x_2,...,x_n]$ be defined by $$h_D(x_1,x_2,...,x_n)=\prod_{\vec{e}=(v,w)\in\vec{E}(D)}(x_v-\varphi_C(e)x_w).$$

We first show that the coefficient of $M_D=\prod_{v\in V}x^{d^+_D(v)}_v$ in $h_D$ is equal to $EE(D_{\sigma,\varphi})-EO(D_{\sigma,\varphi})$. We will want to distinguish between the occurences of a variable $x_v$ when it appears in a factor indexed by $\vec{e}$ such that $v$ is the tail of $\vec{e}$, and the occurences where $v$ is the head of $\vec{e}$. To do so, we introduce a modified polynomial, $\hat{h}_D\in F[x_1,x_2,...,x_n,y_1,y_2,...,y_n]$ where each occurence of a variable $x_v$ in a factor that is indexed by an edge in which $v$ is the tail is replaced by a new variable $y_v$; $$\overline{h}_D(x_1,x_2,...,x_n,y_1,y_2,...,y_n)=\prod_{\vec{e}=(v,w)\in\vec{E}(D)}(y_v-\varphi_C(e)x_w).$$

Furthermore, define $\sigma_C:E(G)\rightarrow \{-1,1\}$ by $\sigma_C(e)=-1$ if $\varphi_C(e)<0$, and $\sigma_C(e)=1$ if $\varphi_C(e)>0$. Then we have 

$$\overline{h}_D(x_1,x_2,...,x_n,y_1,y_2,...,y_n)=\prod_{\vec{e}=(v,w)\in\vec{E}(D)}(y_v-\underbrace{\sigma_C(e)x_w-\sigma_C(e)x_w-\cdots -\sigma_C(e)x_w}_{\varphi^+_C(e) \text{ terms}}).$$

We will also find it useful to denote the $i^{th}$ ocurrence of an $x_w$ term within a $(v,w)$ factor by $x^{(i)}_w$; write 

$$\overline{h}_D(x_1,x_2,...,x_n,y_1,y_2,...,y_n)=\prod_{\vec{e}=(v,w)\in\vec{E}(D)}(y_v-\sigma_C(e)x^{(1)}_w-\sigma_C(e)x^{(2)}_w-\cdots -\sigma_C(e)x^{(\varphi^+_C(e))}_w).$$

Now, consider the monomials in the expansion of $\overline{h}_D$. (Here we are considering the monomials before ``collecting like terms," so there will be $\prod_{e\in E(G)}|1+\varphi^+_C(e)|$ such monomials.)  For each monomial $M$, we can write $M=\prod_{\vec{e}\in \vec{E}}z_{\vec{e}}$, where $z_{\vec{e}}=y_v$ or $z_{\vec{e}}=-\sigma_C(e)x^{(i)}_w$ for some $1\leq i \leq \varphi^+_C(e)$, for each $\vec{e}\in \vec{E}$. 

We will now define a one to one correspondence $M\leftrightarrow S_M$ between the monomials in the expansion of $\overline{h}_D$ and the set of the spanning subdigraphs of $D_{\sigma,\varphi}$ whose edge sets are unions of pairwise edge-disjoint $\gamma$-paths: Given a monomial $M=\prod_{\vec{e}\in \vec{E}}z_{\vec{e}}$ in the expansion of $\overline{h}_D$, as above, associate to $M$ the subdigraph $S_M$ of $D_{\sigma,\varphi}$ defined by $$ V(S_M)=V(D_{\sigma,\varphi})$$
$$\vec{E}(S_M)= \bigcup_{{\vec{e}=(v,w):z_{\vec{e}}=-\sigma_C(e)x^{(i)}_w}{}}\vec{E}(P_{\vec{e},i}) 
.$$

That the above correspondence is indeed a bijection follows from the fact that each monomial is formed by picking exactly one term from each factor of $\overline{h}_D$. Note that the sign of $M$ is negative in the expansion of $\overline{h}_D$ if and only if $$t^-_M=|\{\vec{e}=(v,w)\in \vec{E}(D):z_{\vec{e}}=-\sigma_C(e)x^{(i)}_w, \sigma_C(e)=1\}|$$
 is odd, which occurs if and only the edge set of $S_M$ contains an odd number of $\gamma$-paths of length $3$ (all others have length 4), and thus, $t^-_M$ is odd if and only if $S_M$ has an odd number of edges. 
 
Now, fix  $u\in V(D)$.  Let $M$ be a monomial in the expansion of $\overline{h}_D$ and write $M=\prod_{\vec{e}\in \vec{E}(D)}z_{\vec{e}}$, as above. Define the following notation: 
\begin{center}
\begin{align*}
q_u&=|\{\vec{e}=(u,w)\in \vec{E}(D): z_{\vec{e}}=-\sigma_C(e)x^{(i)}_w, \text{ for some } 1\leq i \leq \varphi^+_C(e) \}|\\
r_u&=|\{\vec{e}=(v,u)\in\vec{E}(D): z_{\vec{e}}=-\sigma_C(e)x^{(i)}_u,  \text{ for some } 1\leq i \leq \varphi^+_C(e)\}|.\\
\end{align*}
\end{center}

Thus, $q_u$ is equal to the number of factors indexed by an edge $\vec{e}=(u,w)$ such that a term of the form $-\sigma_C(e)x^{(i)}_w$ was picked to contribute to $M$, and $r_u$ is the number of factors indexed by an edge $\vec{e}=(v,u)$ such that the a term of the form $-\sigma_C(e)x^{(i)}_u$ was picked to contribute to $M$. By the construction of $S_M$, we have that $d^+_{S_M}(u)=q_u$ and $d^-_{S_M}(u)=r_u$. Hence, $d^+_{S_M}(u)=d^-_{S_M}(u)$ if and only if $ q_u=r_u$. Furthermore, we note that $r_u=\text{deg}_M(x_u)$ (the degree of $x_u$ in the monomial $M$), and since $q_u$ is equal to the number of factors in $\overline{h}_D$ that contain $y_u$ but $-\sigma_C(e)x_w$ was instead picked from that factor to contribute to $M$, and $y_u$ occurs in exactly $d^+_D(u)$ factors of $\overline{h}_D$, we have that $\text{deg}_M(y_u)=d^+_D(u)-q_u$. Therefore, $$d^+_{S_M}(u)=d^-_{S_M}(u) \text{ if and only if } \text{deg}_M(y_u)=d^+_{D}(u)-\text{deg}_M(x_u).$$

This, along with \Cref{prop:gamma}, yields the following characterization of spanning Eulerian subdigraphs of $D_{\sigma,\varphi}$: $S$ is a spanning Eulerian subdigraph of $D_{\sigma,\varphi}$ if and only if $S=S_M$ for some $M$ with $\text{deg}_M(y_v)+\text{deg}_M(x_v)=d^+_{D}(v)$ for every $v\in V(D)$. 

If we now let $y_v=x_v$ in $\overline{h}_D$ for each $v$, we get back the polynomial $h_D$, and the monomials that are in correspondence with the spanning Eulerian subdigraphs of $D_{\sigma,\varphi}$ are now exactly the monomials $M$ in the expansion of $f$ such that for each $v\in V(D)$, $\text{deg}_M(x_v)=d^+_D(v)$; that is, the monomials in correspondence with the spanning Eulerian subdigraphs of $D_{\sigma,\varphi}$ are precisely the occurences of $\pm M_D$ in the expansion of $f$ (recall that $M_D=\prod_{v\in V}x^{d^+_D(v)}_v$). Since letting $y_v=x_v$ does not change the sign of any monomial in the expansion, our earlier observation that $S_M$ has an odd number of edges if and only if $M$ occurs with a negative sign in the expansion of $\overline{h}_D$ tells us that the occurences of $+M_D$ in the expansion of $h_D$ are in correspondence with even spanning Eulerian subdigraphs of $D_{\sigma,\varphi}$ and the occurences of $-M_D$ are in correspondence with odd spanning Eulerian subdigraphs of $D_{\sigma,\varphi}$. Thus, the coefficient of $M_D$ in $h_D$ is $EE(D_{\sigma,\varphi})-EO(D_{\sigma,\varphi})$, which is what we had hoped to show. 

Finally, define $\hat{h}_D\in F[x_1,x_2,...,x_n]$ by $$\hat{h}_D(x_1,x_2,...,x_n)=\prod_{\vec{e}=(v,w)\in\vec{E}(D)}(x_v-\varphi_C(e)x_w-a_{\vec{e}}),$$
so that if $(L,C)$ is an $F$-correspondence assignment on $G$, and $x\in \prod_{i=1}^nL(i)$ with $\hat{h}_D(x)\neq 0$, then $x$ corresponds to an $(L,C)$-coloring of $G$. Then $M_D$ also has maximum degree in $\hat{h}_D$, and the coefficient of $M_D$ in $\hat{h}_D$ is equal to its coefficient in $h_D$. Hence, if $EE(D_{\sigma,\varphi})-EO(D_{\sigma,\varphi})\neq 0$, by the Combinatorial Nullstellensatz (\Cref{thm:CN}), the theorem holds true when $F$ is a field of characteristic 0. 

Now, consider the case where $F=\mathbb{Z}_p$. In this case, as mentioned in the discussion after \Cref{def:zsign}, if $\varphi_C$ is the generalized sign function of $(L,C)$ with respect to $D$, then we can pick any sign function $\sigma_C:E(G)\rightarrow \{-1,1\}$, and pick $\varphi_C^+(e)\in \{1,2,...,p-1\}$ accordingly for each $e\in E(G)$ so that $\varphi_C(e)=\sigma_C(e)\varphi_C^+(e)$, rather than requiring the specific sign function $\sigma_C$ in the above argument for the case where $F$ has characteristic $0$. Once we have done this, we can consider the polynomial $h_D\in  \mathbb{Q}[x_1,x_2,...,x_n]$ as above, and compute its coefficients in $\mathbb{Q}$. Since the coefficients are integers, we can compute the the coefficient of $M_D$ in $\mathbb{Z}_p$ by first computing it in $\mathbb{Q}$ and then reducing it mod $p$.  Thus, the coefficient of $M_D$ in $h_D$, when considering $h_D$ as a polynomial over $\mathbb{Z}_p$ is $EE(D_{\sigma,\varphi})-EO(D_{\sigma,\varphi}) (\text{mod }p)$, and hence, the coefficient of $M_D$ in $\hat{h}_D$ considered as a polynomial over $\mathbb{Z}_p$ is also $EE(D_{\sigma,\varphi})-EO(D_{\sigma,\varphi}) (\text{mod }p)$. The claim therefore follows from the Combinatorial Nullstellensatz in this case, as well.
\end{proof}
\subsection{Orientation theorems for signable and good correspondence assignments}
We can define a simpler auxiliary digraph for signable correspondence assignments and obtain another Alon-Tarsi style theorem as a corollary of \Cref{thm:Eulerianaux}.  

\begin{defi}\label{defi:sigma}
Let $G$ be a multigraph and let $(L,C)$ be a signable $F$-correspondence assigment on $G$, with sign function $\sigma_C:E(G)\rightarrow \{-1,1\}$. For each orientation $D$ of $G$, define $D_{\sigma}$ to be the digraph obtained from $D$ by replacing each edge $\vec{e}=(v,w)\in \vec{E}(D)$ such that $\sigma_C(e)=-1$ with a directed path of length two from $v$ to $w$. 
\end{defi}

 \begin{cor}\label{cor:signaux} Let $G$ be a multigraph, and let $F$ be a field of characteristic $0$ or $F=\mathbb{Z}_p$. Suppose that $(L,C)$ is a signable $F$-correspondence assignment on $G$ with sign function $\sigma_C:E(G)\rightarrow \{-1,1\}$. If there exists an orientation $D$ of $G$ such that $EE(D_\sigma)-EO(D_\sigma)\not\equiv 0_F$ and $|L(v)|\geq d^+_D(v)+1$ for each $v\in V$, then $G$ is $(L,C)$-colorable. 
\end{cor}

\begin{proof} Since $(L,C)$ is signable, we can view $(L,C)$ as a Z-signable correspondence assignment with $\varphi_C:E(G)\rightarrow F\setminus \{0\}$, $\sigma_C:E(G)\rightarrow \{-1,1\}$, and $\varphi^+_C:E\rightarrow F\setminus\{0\}$ such that $\varphi_C=\sigma_C\varphi^+_C $ and $\varphi^+_C(e)=1$ for all $e\in E(G)$. Hence, by \Cref{thm:Eulerianaux}, it suffices to show that in this case $EE(D_\sigma)-EO(D_\sigma)=EE(D_{\sigma,\varphi})-EO(\sigma_{\sigma,\varphi})$. Note that since $\varphi^+_C(e)=1$ for all $e\in E(G)$, in the construction of $D_{\sigma,\varphi}$, an edge $\vec{e}\in \vec{E}(D)$ is replaced by a single directed path of length $4$ when $\sigma_C(e)=-1$, and by a single directed path of length $3$ when $\sigma_C(e)=1$. Since each $\vec{e}\in D$ with $\sigma_C(e)=-1$ corresponds to a directed path of even length in both $D_\sigma$ and $D_{\sigma,\varphi}$, and each $\vec{e}\in \vec{E}(D)$ with $\sigma_C(e)=1$ corresponds to a directed path of odd length in both $D_\sigma$ and $D_{\sigma,\varphi}$, it follows that $EE(D_\sigma)=EE(D_{\sigma,\varphi})$, and $EO(D_\sigma)=EO(D_{\sigma,\varphi})$. Hence  $EE(D_\sigma)-EO(D_\sigma)=EE(D_{\sigma,\varphi})-EO(D_{\sigma,\varphi})$, as desired. 
\end{proof}

\begin{cor}\label{cor:euleriangood} Let $G$ be a multigraph, and let $F$ be a field of characteristic $0$ or $F=\mathbb{Z}_p$. Suppose that $(L,C)$ is a good $F$-correspondence assignment on $G$. Suppose that $D$ is an orientation of $G$ such that $EE(D)-EO(D)\not\equiv 0_F$. If $|L(v)|\geq d^+_D(v)+1$ for each $v\in V$, then $G$ is $(L,C)$-colorable. 
\end{cor}
\begin{proof}
Since $(L,C)$ is a good correspondence assignment, we can think of $(L,C)$ as a signable correspondence assignment on $G$, with $\sigma(e)=1$ for all $e\in E(G)$. By the construction of $D_{\sigma}$, we then have that $D_\sigma=D$. Hence, the claim follows from \Cref{cor:signaux}.

\end{proof}
\begin{cor} \label{cor:subdivide}
Let $G$ be a multigraph, $F$ a field of characteristic $0$, and $(L,C)$ a signable $F$-correspondence assignment on $G$. Suppose that subdividing each bad edge of $G$ results in a bipartite graph. If there exists an orientation $D$ of $G$ such that $|L(v)|\geq d^+_D(v)+1$ for all $v\in V$, Then $G$ is $(L,C)$-colorable.  
\end{cor}

\begin{proof}
Since subdividing every bad edge yields a bipartite graph, we have that $EE(D_\sigma)>0$ and $EO(D_\sigma)=0$. Hence, by \Cref{cor:signaux}, the claim holds true.
\end{proof}

\begin{exam} As an illustration of \Cref{cor:subdivide}, let us consider the simple example of the the wheel graph, $W_6$. Note that the orientation of $W_6$ and the list assignment shown below demonstrate that it is not the case that for an arbitrary orientation $D$ of $W_6$, if $(L,C)$ is an $\mathbb{R}$-correspondence assignment on $W_6$ with $|L(v)|\geq d^+_D(v)+1$, then $W_6$ is $(L,C)$-colorable, since $W_6$ is not colorable from the lists.  

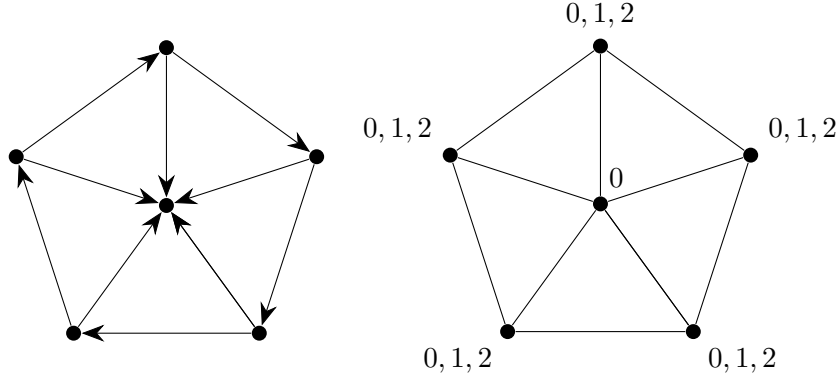
\begin{center}
\begin{tikzpicture}[scale=.7, >={Stealth[scale=1.8]}, 
dot/.style={circle,fill,inner sep=2pt},
    declare function={R=3;},bend angle=12]
 \path[dash pattern=on 1.5pt off 1pt] 
 foreach \X [count=\Y] in {1,2,3,4,5}
 {(162-\Y*72:R) node[dot,label={162-\Y*72:
 }] 
 (v\X){}
 \ifnum\Y>1
  \ifnum\Y<6
   (v\the\numexpr\Y-1) edge[solid,->] (v\Y)
  \fi
 \fi
 };
 \path[->] 
(v5) edge (v1);

\node[] at (0,-3.17){};
 \node[circle,fill,inner sep=2pt] (M) at (0,0){};
\path[<-]
(M) edge (v5)

(M) edge (v1)
(M) edge (v3)
;
 \path[<-]
(M) edge (v4)
(M) edge (v3)
(M) edge (v2)
;
\end{tikzpicture}
\begin{tikzpicture}[scale=.7, 
dot/.style={circle,fill,inner sep=2pt},
    declare function={R=3;},bend angle=12]
 \path[dash pattern=on 1.5pt off 1pt] 
 foreach \X [count=\Y] in {1,2,3,4,5}
 {(162-\Y*72:R) node[dot,label={162-\Y*72:$0,1,2$ 
 }] 
 (v\X){}
 \ifnum\Y>1
  \ifnum\Y<6
   (v\the\numexpr\Y-1) edge[solid,-] (v\Y)
  \fi
 \fi
 };
 \path[-] 
(v5) edge (v1);

 \node[circle,fill,inner sep=2pt] (M) at (0,0){};
 \node[] at (0.3,0.5){$0$};
\path[-]
(M) edge (v5)

(M) edge (v1)
(M) edge (v3)
;
 \path[-]
(M) edge (v4)
(M) edge (v3)
(M) edge (v2)
;
\end{tikzpicture}
\end{center}

However, suppose that $(L,C)$ is a signable $\mathbb{R}$-correspondence assignment such that each solid edge in the diagram below is a good edge, and each dashed edge is a bad edge. Then subdividing the dashed edges yields a bipartite graph, so if there exists any orientation $D$ of $W_6$ such that $|L(v)|\geq d^+_D(v)+1$ for each $v\in V(W_6)$, then $W_6$ is $(L,C)$-colorable by \Cref{cor:subdivide}. 

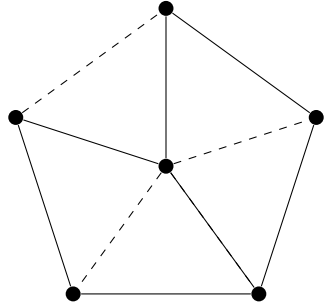
\begin{center}
\begin{tikzpicture}[scale=.7, 
dot/.style={circle,fill,inner sep=2pt},
    declare function={R=3;},bend angle=12]
 \path[dash pattern=on 1.5pt off 1pt] 
 foreach \X [count=\Y] in {1,2,3,4,5}
 {(162-\Y*72:R) node[dot,label={162-\Y*72:
 }] 
 (v\X){}
 \ifnum\Y>1
  \ifnum\Y<6
   (v\the\numexpr\Y-1) edge[solid,-] (v\Y)
  \fi
 \fi
 };
 \path[-,dashed] 
(v5) edge (v1);

 \node[circle,fill,inner sep=2pt] (M) at (0,0){};
\path[-]
(M) edge (v5)

(M) edge (v1)
(M) edge (v3)
;
 \path[-, dashed]
(M) edge (v4)
(M) edge (v3)
(M) edge (v2);
\end{tikzpicture}

\end{center}

\end{exam}

\section{Orientation theorems for DP-Colorings of simple graphs}

In this section, we use the results from Section 4 to prove some Alon-Tarsi style theorems for DP-colorings of simple graphs. The main idea behind these theorems is that for an arbitrary $F$-correspondence $(L,C)$ on a graph $G$, we can always construct a suitable multigraph $G'$, and a correspondence assignment $(L,C')$ on $G'$, that is either good, signable, or Z-signable, and so that $G$ is $(L,C)$-colorable if and only if $G'$ is $(L,C')$-colorable. We can then apply the theorems proved in Section 4 to the multigraph $G'$.

Let $(L,C)$ be an $F$-correspondence assignment on a graph $G$ and fix an orientation $D$ of $G$. It will be useful to apply our terminology from the definitions in Section 3 to describe the matchings $C_e$ as follows: Let $e=\{v,w\}\in E(G)$ such that $\vec{e}=(v,w)\in \vec{E}(D)$.  The matching $C_e$ is \textbf{good} if there exists some $a_{\vec{e}}\in F$  such that  $c_1-c_2=a_{\vec{e}}$ whenever $\{(v,c_1),(w,c_2)\}\in C_e$. The matching $C_e$ is \textbf{signable} if there exists some $\sigma(e)\in\{-1,1\}$ and some $a_{\vec{e}}\in F$  such that $c_1-\sigma(e)c_2=a_{\vec{e}}$ whenever $\{(v,c_1),(w,c_2)\}\in C_e$. The matching $C_{\{v,w\}}$ is \textbf{Z-signable} with respect to the orientation $v \rightarrow w$ of $\{v,w\}$, if there exists some $\varphi(e)\in\langle 1 \rangle \setminus \{0\}$ and some $a_{\vec{e}}\in F$  such that $c_1-\varphi(e)c_2=a_{\vec{e}}$ for all $\{(v,c_1),(w,c_2)\}\in C_{e}$. (Recall that $\langle 1 \rangle$ is the additive subgroup of $F$ generated by $1$.)

\begin{defi}{\label{def:minz}}
Suppose that $G=(V,E)$ is a simple graph and let $(L,C)$ be an $F$-correspondence assignment on $G$. For each $e\in E(G)$, define $\omega^z(e)$ to be the least $k$ such that $C_e$ is the union of $k$ partial matchings, $C^z_{e_1},C^z_{e_2},...,C^z_{e_k}$, each of which is  $Z$-signable. Define $G^z$ to be the multigraph formed by replacing each edge $e\in E(G)$ with $\omega^z(e)$ edges $e_1,...,e_{\omega^z(e)}$.
\end{defi}

Let $(L,C)$ be an $F$-correspondence on a simple graph $G$. For each $e\in E(G)$, fix a set of $\omega^z(e)$ matchings, $C^z_{e_1},C^z_{e_2},...,C^z_{e_{\omega^z(e)}}$, such that $C_e=\bigcup_{i=1}^{\omega^Z(e)}C^z_{e_i}$, and  $C^z_{e_i}$ is $Z$-signable for each $i$. Define the $Z$-signable $F$-correspondence assignment $(L,C^z)$ on $G^z$ so that the matchings of $C^z$ are given by $C^z_{e_{i}}$, for each $e_{i}\in E(G^z)$.

\begin{cor}\label{cor:multizsign}
Let $G=(V,E)$ be a simple graph and let $F$ be a field of characteristic $0$ or $F=\mathbb{Z}_p$. Let $(L,C)$ be an $F$-correspondence assignment on $G$. Let $G^z$ and $(L,C^z)$ be defined as above. If there exists an orientation $D$ of $G^z$ such that each edge of  $G^z$ is $Z$-signable with respect to $D$, $EE(D_{\sigma,\varphi})-EO(D_{\sigma,\varphi})\not \equiv 0_F$, and such that $|L(v)|\geq d^+_D(v)+1$ for each $v\in V(G)$, then $G$ is $(L,C)$ colorable. 
\end{cor}

\begin{proof}
    Since for each $e\in E(G)$, $C_e=\bigcup_{i=1}^{\omega^z(e)}C^z_{e_i}$, we have that $G$ is $(L,C)$-colorable if and only if $G^z$ is $(L,C^z)$-colorable. By \Cref{thm:Eulerianaux}, if $EE(D_{\sigma,\varphi})-EO(D_{\sigma,\varphi})\not \equiv 0_F$, then $G^z$ is $(L,C^z)$-colorable. Hence $G$ is $(L,C)$-colorable, as desired. 
\end{proof}

\begin{defi}{\label{def:mins}}
Suppose that $G$ is a simple graph, let $F$ be a field of characteristic $0$ or $F=\mathbb{Z}_p$, and let $(L,C)$ be an $F$-correspondence assignment on $G$. For each $e\in E(G)$, define $\omega^s(e)$ to be the least $k$ such that $C_e$ is the union of $k$ partial matchings, $C^s_{e_1},C^s_{e_2},...,C^s_{e_k}$, each of which is signable. Define $G^s$ to be the multigraph formed by replacing each edge $e\in E(G)$ with $\omega^s(e)$ edges $e_1,...,e_{\omega^s(e)}$.
\end{defi}
 
 Analogously to the $Z$-signable case, given an $F$-correspondence assignment $(L,C)$ on $G$, for each $e\in E(G)$, fix a set of $\omega^s(e)$ matchings, $C^s_{e_1},C^s_{e_2},...,C^s_{e_{\omega^s(e)}}$, such that $C_e=\bigcup_{i=1}^{\omega^s(e)}C^s_{e_i}$, and  $C^s_{e_i}$ is signable for each $i$. Define the signable $F$-correspondence assignment $(L,C^s)$ on $G^s$ so that the matchings of $C^s$ are given by $C^s_{e_{i}}$, for each $e_{i}\in E(G^s)$. 

 \begin{cor}\label{cor:multisign}
 Let $G=(V,E)$ be a simple graph and let let $F$ be a field of characteristic $0$ or $F=\mathbb{Z}_p$. Let $(L,C)$ be an  $F$-correspondence assignment on $G$. Let $G^s$ and $(L,C^s)$ be defined as above. If there exists an orientation $D$ of $G^s$ such that $EE(D_{\sigma})-EO(D_{\sigma})\not \equiv 0_F$, and such that $|L(v)|\geq d^+_D(v)+1$ for each $v\in V(G)$, then $G$ is $(L,C)$ colorable. 
\end{cor}

\begin{proof}
As in the proof of \Cref{cor:signaux}, since $(L,C^s)$ is signable, we can view $(L,C^s)$ as a Z-signable correspondence assignment with $\varphi_C:E(G)\rightarrow F\setminus \{0\}$, $\sigma_C:E(G)\rightarrow \{-1,1\}$, and $\varphi^+_C:E\rightarrow F\setminus\{0\}$ such that $\varphi_C=\sigma_C\varphi^+_C $ and $\varphi^+_C(e)=1$ for all $e\in E(G)$. We showed in \Cref{cor:signaux} that in this case, $EE(D_\sigma)-EO(D_\sigma)=EE(D_{\sigma,\varphi})-EO(\sigma_{\sigma,\varphi})$.  Therefore, the claim follows from \Cref{cor:multizsign}.
\end{proof}

\begin{defi}\label{def:ming}Let $G=(V,E)$, with $V=\{1,2,..,n\}$, let $F$ be a field, and let $(L,C)$ be an $F$-correspondence assignment on $G$. For $e=\{i,j\}\in E(G)$, with $i<j$, define $\omega^g(e)$ by $$\omega^g(e)=|\{a \in F:\text{ there exists some } \{(i,c_1),(j,c_2)\}\in C_{\{i,j\}}\text{ such that }a=c_1-c_2\}|.$$
 Define $G^g$ to be the multigraph formed by replacing each edge $e\in E(G)$ with $\omega^g(e)$ edges $e_1,...,e_{\omega^g(e)}$
\end{defi}

Note that $\omega^g(e)$ is equal to the least $k$ such that $C_e$ can be written as a union of $k$ good matchings. Similarly to the $Z$-signable and signable cases above, given an $F$-correspondence assignment $(L,C)$ on $G$, for each $e\in E(G)$, fix a set of $\omega^g(e)$ matchings, $C^g_{e_1},C^g_{e_2},...,C^g_{e_{\omega^g(e)}}$, such that $C_e=\bigcup_{i=1}^{\omega^g(e)}C^g_{e_i}$, and  $C^g_{e_i}$ is good for each $i$. Define the good $F$-correspondence assignment $(L,C^g)$ on $G^g$ so that the matchings of $C^g$ are given by $C^g_{e_{i}}$, for each $e_{i}\in E(G^g)$.

 \begin{cor}\label{cor:multigood}
 Let $G=(V,E)$ be a simple graph and let let $F$ be a field of characteristic $0$ or $F=\mathbb{Z}_p$. Let $(L,C)$ be an $F$-correspondence assignment on $G$. Let $G^g$ and $(L,C^g)$ be defined as above. If there exists an orientation $D$ of $G^g$ such that $EE(D)-EO(D)\not \equiv 0_F$, and such that $|L(v)|\geq d^+_D(v)+1$ for each $v\in V(G)$, then $G$ is $(L,C)$-colorable. 
\end{cor}

\begin{proof}
Since $(L,C^g)$ is a good correspondence assignment on $G^g$, if there exists an orientation $D$ of $G$ such that $EE(D)-EO(D)\not \equiv 0_F$, then then $G^g$ is $(L,C^g)$-colorable, by \Cref{cor:euleriangood}. 
Then, since $G^g$ is $(L,C^g)$-colorable if and only if $G$ is $(L,C)$-colorable, we have that if $G$ is $(L,C)$-colorable, as desired.  
\end{proof}

\begin{cor}\label{cor:absval}
Let $G$ be a simple graph with $V(G)=\{1,2,...,n\}$, and let $(L,C)$ be an $\mathbb{R}$-correspondence on $G$ such that for each edge $\{i,j\}\in E(G)$ with $i<j$ that is not a good edge, there exists some $a_{ij} \in \mathbb{R}$ such that  $|c_1-c_2|=a_{ij}$ whenever $\{(i,c_1),(j,c_2)\}\in C_{\{i,j\}}$. Let $G'$ be the multigraph constructed from $G$ by replacing each edge of $G$ that is not good by two edges. If there exists an orientation $D$ of $G$ with $EE(D)\neq EO(D)$ and $|L(v)|\geq d^+_D(v)+1$ for each $v\in V$, then $G$ is $(L,C)$ colorable.
\end{cor}

\begin{proof} If $e\in E(G)$ is not a good edge, we have $\omega^g(e)=2$, and if $e$ is a good edge, we have and $\omega^g(e)=1$. Thus, the corollary follows from \Cref{cor:multigood}. 

\end{proof}

Note that \Cref{cor:absval} applies to the example in \Cref{fig:badcor}. We give another example that illustrates how \Cref{cor:absval} might be applied.

\begin{exam}\label{exam:toroidalgrid}
Consider the following toroidal grid, $G$. (Here, the torus is drawn as a square with opposite sides identified with one another and the vertices and edges on the outer boundary are identified with the corresponding vertices and edges on the sides opposite them.) Suppose that $(L,C)$ is an $\mathbb{R}$-correspondence on $G$ such that every solid edge in the diagram below is good and for every dashed edge $e=\{i,j\}$, $i<j$ there exists some $a_{ij}\in \mathbb{R}$ such that $\{(i,c_1),(j,c_2)\}\in C_{e}$ implies $|c_1- c_2|=a_{ij}$. Then $G^g$ is shown below. Since $G^g$ is bipartite, $EE(D)>0$ and $EO(D)=0$ for any orientation $D$ of $G^g$, so by \Cref{cor:absval}, if there exists some orientation $D$ of $G^g$ such that $|L(v)|\geq d^+_{D}(v)+1$ for each $v\in V$, then we know that $G$ is $(L,C)$ colorable.  
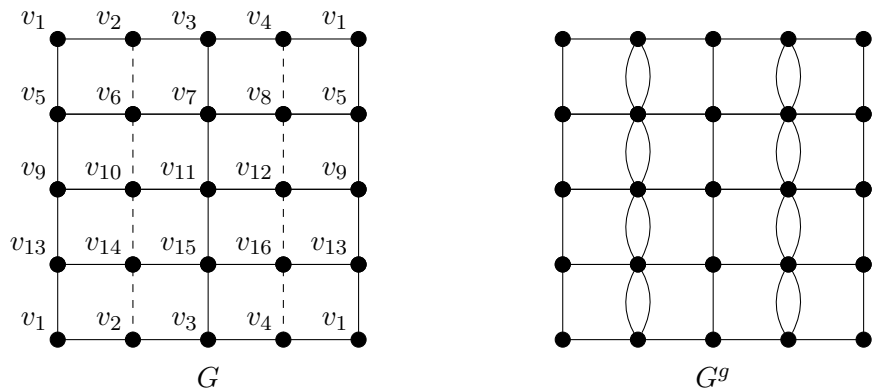
\vspace{.2cm}

\begin{center}\begin{tikzpicture}[scale=1]

\tikzset{vertex/.style = {shape=circle,fill=black, draw, inner sep=0pt,minimum size=2 mm}}
\tikzset{edge/.style = {-
}}
\node[] at (2,-.5){${G}$};

\node[anchor=south east] at  (0,0) {$v_1$};
\node[anchor=south east] at  (1,0) {$v_2$};
\node[anchor=south east] at  (2,0) {$v_3$};
\node[anchor=south east] at  (3,0) {$v_4$};
\node[anchor=south east] at  (4,0) {$v_1$};

\node[anchor=south east] at  (0,4) {$v_1$};
\node[anchor=south east] at  (1,4) {$v_2$};
\node[anchor=south east] at  (2,4) {$v_3$};
\node[anchor=south east] at  (3,4) {$v_4$};
\node[anchor=south east] at  (4,4) {$v_1$};

\node[anchor=south east] at  (0,3) {$v_5$};
\node[anchor=south east] at  (1,3) {$v_6$};
\node[anchor=south east] at  (2,3) {$v_7$};
\node[anchor=south east] at  (3,3) {$v_8$};
\node[anchor=south east] at  (4,3) {$v_5$};


\node[anchor=south east] at  (0,2) {$v_9$};
\node[anchor=south east] at  (1,2) {$v_{10}$};
\node[anchor=south east] at  (2,2) {$v_{11}$};
\node[anchor=south east] at  (3,2) {$v_{12}$};
\node[anchor=south east] at  (4,2) {$v_9$};

\node[anchor=south east] at  (0,1) {$v_{13}$};
\node[anchor=south east] at  (1,1) {$v_{14}$};
\node[anchor=south east] at  (2,1) {$v_{15}$};
\node[anchor=south east] at  (3,1) {$v_{16}$};
\node[anchor=south east] at  (4,1) {$v_{13}$};


\node[vertex] (0,0) at  (0,0) {};
\node[vertex] (0,1) at  (0,1) {};
\node[vertex] (1,1) at  (1,1) {};
\node[vertex] (1,0) at  (1,0) {};

\draw[edge] (0,0) to (1,0);
\draw[edge] (0,1) to (0,0);
\draw[edge, dashed] (1,0) to (1,1);

\draw[edge] (0,1) to (1,1);


\node[vertex] (2,1) at  (2,1) {};
\node[vertex] (2,0) at  (2,0) {};

\draw[edge] (1,0) to (2,0);
\draw[edge] (2,1) to (2,0);
\draw[edge] (1,1) to (2,1);

\node[vertex] (3,1) at  (3,1) {};
\node[vertex] (3,0) at  (3,0) {};

\draw[edge] (2,0) to (3,0);
\draw[edge, dashed] (3,0) to (3,1);
\draw[edge] (2,1) to (3,1);

\node[vertex] (4,1) at  (4,1) {};
\node[vertex] (4,0) at  (4,0) {};

\draw[edge] (3,0) to (4,0);
\draw[edge] (4,1) to (4,0);
\draw[edge] (3,1) to (4,1);


\node[vertex] (0,1) at  (0,1) {};
\node[vertex] (0,2) at  (0,2) {};
\node[vertex] (1,2) at  (1,2) {};
\node[vertex] (1,1) at  (1,1) {};
\draw[edge] (0,1) to (1,1);
\draw[edge] (0,2) to (0,1);
\draw[edge,dashed ] (1,1) to (1,2);
\draw[edge] (0,2) to (1,2);


\node[vertex] (2,2) at  (2,2) {};
\node[vertex] (2,1) at  (2,1) {};

\draw[edge] (1,1) to (2,1);
\draw[edge] (2,2) to (2,1);
\draw[edge] (1,2) to (2,2);

\node[vertex] (3,2) at  (3,2) {};
\node[vertex] (3,1) at  (3,1) {};

\draw[edge] (2,1) to (3,1);
\draw[edge,dashed] (3,1) to (3,2);
\draw[edge] (2,2) to (3,2);

\node[vertex] (4,2) at  (4,2) {};
\node[vertex] (4,1) at  (4,1) {};

\draw[edge] (3,1) to (4,1);
\draw[edge] (4,2) to (4,1);
\draw[edge] (3,2) to (4,2);


\node[vertex] (0,2) at  (0,2) {};
\node[vertex] (0,3) at  (0,3) {};
\node[vertex] (1,3) at  (1,3) {};
\node[vertex] (1,2) at  (1,2) {};
\draw[edge] (0,2) to (1,2);
\draw[edge] (0,3) to (0,2);
\draw[edge,dashed] (1,2) to (1,3);
\draw[edge] (0,3) to (1,3);


\node[vertex] (2,3) at  (2,3) {};
\node[vertex] (2,2) at  (2,2) {};

\draw[edge] (1,2) to (2,2);
\draw[edge] (2,3) to (2,2);
\draw[edge] (1,3) to (2,3);

\node[vertex] (3,3) at  (3,3) {};
\node[vertex] (3,2) at  (3,2) {};

\draw[edge] (2,2) to (3,2);
\draw[edge,dashed] (3,2) to (3,3);
\draw[edge] (2,3) to (3,3);

\node[vertex] (4,3) at  (4,3) {};
\node[vertex] (4,2) at  (4,2) {};

\draw[edge] (3,2) to (4,2);
\draw[edge] (4,3) to (4,2);
\draw[edge] (3,3) to (4,3);

\node[vertex] (0,3) at  (0,3) {};
\node[vertex] (0,4) at  (0,4) {};
\node[vertex] (1,4) at  (1,4) {};
\node[vertex] (1,3) at  (1,3) {};
\draw[edge] (0,3) to (1,3);
\draw[edge] (0,4) to (0,3);
\draw[edge,dashed] (1,3) to (1,4);
\draw[edge] (0,4) to (1,4);


\node[vertex] (2,4) at  (2,4) {};
\node[vertex] (2,3) at  (2,3) {};

\draw[edge] (1,3) to (2,3);
\draw[edge] (2,4) to (2,3);
\draw[edge] (1,4) to (2,4);

\node[vertex] (3,4) at  (3,4) {};
\node[vertex] (3,3) at  (3,3) {};

\draw[edge] (2,3) to (3,3);
\draw[edge,dashed] (3,3) to (3,4);
\draw[edge] (2,4) to (3,4);

\node[vertex] (4,4) at  (4,4) {};
\node[vertex] (4,3) at  (4,3) {};

\draw[edge] (3,3) to (4,3);
\draw[edge] (4,4) to (4,3);
\draw[edge] (3,4) to (4,4);
\end{tikzpicture}\hspace{2.5cm}
\begin{tikzpicture}[scale=1]

\tikzset{vertex/.style = {shape=circle,fill=black, draw, inner sep=0pt,minimum size=2 mm}}
\tikzset{edge/.style = {-
}}

\node[] at (2,-.5){${G^g}$};

\node[vertex] (0,0) at  (0,0) {};
\node[vertex] (0,1) at  (0,1) {};
\node[vertex] (1,1) at  (1,1) {};
\node[vertex] (1,0) at  (1,0) {};
\draw[edge] (0,0) to (1,0);
\draw[edge] (0,1) to (0,0);
\draw[edge] (1,0) to[bend left] (1,1);
\draw[edge] (1,1) to[bend left] (1,0);
\draw[edge] (0,1) to (1,1);


\node[vertex] (2,1) at  (2,1) {};
\node[vertex] (2,0) at  (2,0) {};

\draw[edge] (1,0) to (2,0);
\draw[edge] (2,1) to (2,0);
\draw[edge] (1,1) to (2,1);

\node[vertex] (3,1) at  (3,1) {};
\node[vertex] (3,0) at  (3,0) {};

\draw[edge] (2,0) to (3,0);
\draw[edge] (3,0) to[bend left] (3,1);
\draw[edge] (3,1) to[bend left] (3,0);
\draw[edge] (2,1) to (3,1);

\node[vertex] (4,1) at  (4,1) {};
\node[vertex] (4,0) at  (4,0) {};

\draw[edge] (3,0) to (4,0);
\draw[edge] (4,1) to (4,0);
\draw[edge] (3,1) to (4,1);


\node[vertex] (0,1) at  (0,1) {};
\node[vertex] (0,2) at  (0,2) {};
\node[vertex] (1,2) at  (1,2) {};
\node[vertex] (1,1) at  (1,1) {};
\draw[edge] (0,1) to (1,1);
\draw[edge] (0,2) to (0,1);
\draw[edge] (1,1) to[bend left] (1,2);
\draw[edge] (1,2) to[bend left] (1,1);
\draw[edge] (0,2) to (1,2);


\node[vertex] (2,2) at  (2,2) {};
\node[vertex] (2,1) at  (2,1) {};

\draw[edge] (1,1) to (2,1);
\draw[edge] (2,2) to (2,1);
\draw[edge] (1,2) to (2,2);

\node[vertex] (3,2) at  (3,2) {};
\node[vertex] (3,1) at  (3,1) {};

\draw[edge] (2,1) to (3,1);
\draw[edge] (3,1) to[bend left] (3,2);
\draw[edge] (3,2) to[bend left] (3,1);
\draw[edge] (2,2) to (3,2);

\node[vertex] (4,2) at  (4,2) {};
\node[vertex] (4,1) at  (4,1) {};

\draw[edge] (3,1) to (4,1);
\draw[edge] (4,2) to (4,1);
\draw[edge] (3,2) to (4,2);


\node[vertex] (0,2) at  (0,2) {};
\node[vertex] (0,3) at  (0,3) {};
\node[vertex] (1,3) at  (1,3) {};
\node[vertex] (1,2) at  (1,2) {};
\draw[edge] (0,2) to (1,2);
\draw[edge] (0,3) to (0,2);
\draw[edge] (1,2) to[bend left] (1,3);
\draw[edge] (1,3) to[bend left] (1,2);
\draw[edge] (0,3) to (1,3);


\node[vertex] (2,3) at  (2,3) {};
\node[vertex] (2,2) at  (2,2) {};

\draw[edge] (1,2) to (2,2);
\draw[edge] (2,3) to (2,2);
\draw[edge] (1,3) to (2,3);

\node[vertex] (3,3) at  (3,3) {};
\node[vertex] (3,2) at  (3,2) {};

\draw[edge] (2,2) to (3,2);
\draw[edge] (3,2) to[bend left] (3,3);
\draw[edge] (3,3) to[bend left] (3,2);
\draw[edge] (2,3) to (3,3);

\node[vertex] (4,3) at  (4,3) {};
\node[vertex] (4,2) at  (4,2) {};

\draw[edge] (3,2) to (4,2);
\draw[edge] (4,3) to (4,2);
\draw[edge] (3,3) to (4,3);

\node[vertex] (0,3) at  (0,3) {};
\node[vertex] (0,4) at  (0,4) {};
\node[vertex] (1,4) at  (1,4) {};
\node[vertex] (1,3) at  (1,3) {};
\draw[edge] (0,3) to (1,3);
\draw[edge] (0,4) to (0,3);
\draw[edge] (1,3) to[bend left] (1,4);
\draw[edge] (1,4) to[bend left] (1,3);
\draw[edge] (0,4) to (1,4);


\node[vertex] (2,4) at  (2,4) {};
\node[vertex] (2,3) at  (2,3) {};

\draw[edge] (1,3) to (2,3);
\draw[edge] (2,4) to (2,3);
\draw[edge] (1,4) to (2,4);

\node[vertex] (3,4) at  (3,4) {};
\node[vertex] (3,3) at  (3,3) {};

\draw[edge] (2,3) to (3,3);
\draw[edge] (3,3) to[bend left] (3,4);
\draw[edge] (3,4) to[bend left] (3,3);
\draw[edge] (2,4) to (3,4);

\node[vertex] (4,4) at  (4,4) {};
\node[vertex] (4,3) at  (4,3) {};

\draw[edge] (3,3) to (4,3);
\draw[edge] (4,4) to (4,3);
\draw[edge] (3,4) to (4,4);
\end{tikzpicture}
\end{center} 

In particular, the existence of the following orientation $D'$ of $G^g$ implies that if $(L,C)$ is an $\mathbb{R}$-correspondence satisfying the required conditions, and $|L(v_i)|\geq 4$ whenever $i$ is even, and $|L(v_i)|\geq 3$ whenever $i$ is odd, then 
$G$ is $(L,C)$ colorable.
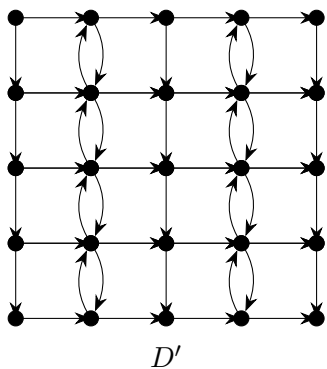
\begin{center}
    \begin{tikzpicture}[scale=1]

\tikzset{vertex/.style = {shape=circle,fill=black, draw, inner sep=0pt,minimum size=2 mm}}
\tikzset{edge/.style = {-{Stealth[scale=1.3]}
}}

\node[] at (2,-.5){${D'}$};

\node[vertex] (0,0) at  (0,0) {};
\node[vertex] (0,1) at  (0,1) {};
\node[vertex] (1,1) at  (1,1) {};
\node[vertex] (1,0) at  (1,0) {};
\draw[edge] (0,0) to (1,0);
\draw[edge] (0,1) to (0,0);
\draw[edge] (1,0) to[bend left] (1,1);
\draw[edge] (1,1) to[bend left] (1,0);
\draw[edge] (0,1) to (1,1);


\node[vertex] (2,1) at  (2,1) {};
\node[vertex] (2,0) at  (2,0) {};

\draw[edge] (1,0) to (2,0);
\draw[edge] (2,1) to (2,0);
\draw[edge] (1,1) to (2,1);

\node[vertex] (3,1) at  (3,1) {};
\node[vertex] (3,0) at  (3,0) {};

\draw[edge] (2,0) to (3,0);
\draw[edge] (3,0) to[bend left] (3,1);
\draw[edge] (3,1) to[bend left] (3,0);
\draw[edge] (2,1) to (3,1);

\node[vertex] (4,1) at  (4,1) {};
\node[vertex] (4,0) at  (4,0) {};

\draw[edge] (3,0) to (4,0);
\draw[edge] (4,1) to (4,0);
\draw[edge] (3,1) to (4,1);


\node[vertex] (0,1) at  (0,1) {};
\node[vertex] (0,2) at  (0,2) {};
\node[vertex] (1,2) at  (1,2) {};
\node[vertex] (1,1) at  (1,1) {};
\draw[edge] (0,1) to (1,1);
\draw[edge] (0,2) to (0,1);
\draw[edge] (1,1) to[bend left] (1,2);
\draw[edge] (1,2) to[bend left] (1,1);
\draw[edge] (0,2) to (1,2);


\node[vertex] (2,2) at  (2,2) {};
\node[vertex] (2,1) at  (2,1) {};

\draw[edge] (1,1) to (2,1);
\draw[edge] (2,2) to (2,1);
\draw[edge] (1,2) to (2,2);

\node[vertex] (3,2) at  (3,2) {};
\node[vertex] (3,1) at  (3,1) {};

\draw[edge] (2,1) to (3,1);
\draw[edge] (3,1) to[bend left] (3,2);
\draw[edge] (3,2) to[bend left] (3,1);
\draw[edge] (2,2) to (3,2);

\node[vertex] (4,2) at  (4,2) {};
\node[vertex] (4,1) at  (4,1) {};

\draw[edge] (3,1) to (4,1);
\draw[edge] (4,2) to (4,1);
\draw[edge] (3,2) to (4,2);


\node[vertex] (0,2) at  (0,2) {};
\node[vertex] (0,3) at  (0,3) {};
\node[vertex] (1,3) at  (1,3) {};
\node[vertex] (1,2) at  (1,2) {};
\draw[edge] (0,2) to (1,2);
\draw[edge] (0,3) to (0,2);
\draw[edge] (1,2) to[bend left] (1,3);
\draw[edge] (1,3) to[bend left] (1,2);
\draw[edge] (0,3) to (1,3);


\node[vertex] (2,3) at  (2,3) {};
\node[vertex] (2,2) at  (2,2) {};

\draw[edge] (1,2) to (2,2);
\draw[edge] (2,3) to (2,2);
\draw[edge] (1,3) to (2,3);

\node[vertex] (3,3) at  (3,3) {};
\node[vertex] (3,2) at  (3,2) {};

\draw[edge] (2,2) to (3,2);
\draw[edge] (3,2) to[bend left] (3,3);
\draw[edge] (3,3) to[bend left] (3,2);
\draw[edge] (2,3) to (3,3);

\node[vertex] (4,3) at  (4,3) {};
\node[vertex] (4,2) at  (4,2) {};

\draw[edge] (3,2) to (4,2);
\draw[edge] (4,3) to (4,2);
\draw[edge] (3,3) to (4,3);

\node[vertex] (0,3) at  (0,3) {};
\node[vertex] (0,4) at  (0,4) {};
\node[vertex] (1,4) at  (1,4) {};
\node[vertex] (1,3) at  (1,3) {};
\draw[edge] (0,3) to (1,3);
\draw[edge] (0,4) to (0,3);
\draw[edge] (1,3) to[bend left] (1,4);
\draw[edge] (1,4) to[bend left] (1,3);
\draw[edge] (0,4) to (1,4);


\node[vertex] (2,4) at  (2,4) {};
\node[vertex] (2,3) at  (2,3) {};

\draw[edge] (1,3) to (2,3);
\draw[edge] (2,4) to (2,3);
\draw[edge] (1,4) to (2,4);

\node[vertex] (3,4) at  (3,4) {};
\node[vertex] (3,3) at  (3,3) {};

\draw[edge] (2,3) to (3,3);
\draw[edge] (3,3) to[bend left] (3,4);
\draw[edge] (3,4) to[bend left] (3,3);
\draw[edge] (2,4) to (3,4);

\node[vertex] (4,4) at  (4,4) {};
\node[vertex] (4,3) at  (4,3) {};

\draw[edge] (3,3) to (4,3);
\draw[edge] (4,4) to (4,3);
\draw[edge] (3,4) to (4,4);
\end{tikzpicture}

\end{center}

\end{exam} 

Note that the above example can be easily generalized to $2k\times 2k$ toroidal grids, rather than just the $4\times 4$ grid. We further note that the Alon-Tarsi method for list colorings has been applied recently to toroidal grids and cartesian products by others, for example in  \cite{MR3904819} and \cite{MR4602837}. We wonder if these results might be extended to broader classes of correspondence assignments such as those mentioned in this work.

\section{Concluding remarks}
We first remark on limitations of the results in Section 5. Though the theorems in Section $5$ do indeed give sufficient combinatorial conditions for a graph $G$ to be $(L,C)$-colorable, if the matchings $C_e$ are not well behaved (in the sense that they are not mostly generalized signable), the bounds on the list length given by the out-degree of the orientations quickly become worse than known bounds. In particular, the results on degree-list DP-colorability in \cite{MR3686937} are often better than what is produced by the orientation theorems in Section 5. We remark, however, that for correspondence assignments such as those given in \Cref{exam:toroidalgrid}, we do not know of an easier way to prove (L,C)-colorability. 

 Another point of consideration is that, though the results in Sections 3,4, and 5 all concern the $(L,C)$ colorability of a graph $G$ for $F$-correspondence assignments $(L,C)$, any result that proves $(L,C)$-colorability of $G$ also proves $(L',C')$-colorability of $G$ for any equivalent correspondence assignment $(L',C')$, even though $(L',C')$ may not satisfy the hypothesis of a given theorem before renaming its list elements. Furthermore, the list elements of $(L',C')$ may come from any set, not just a field $F$. Therefore, it would be interesting to identify classes of graphs and correspondence assignments that allow for methodical renamings of their list elements, so that we obtain an $F$-correspondence assignment with desireable field theoretic properties, as in the classes introduced in this paper.   

Finally, we mention a possible practical use for the theorems in this work.  As mentioned in the introduction, one of the overarching goals in the study of DP-colorings is to determine similarities and differences between DP-colorings and list colorings. Results that highlight differences between the two often require finding only one example of a correspondence assignment from which a graph is not DP-colorable. Thus, to help guide our efforts, it is useful to have theorems that characterize classes of correspondence assignments that must admit a coloring, as is the case with the theorems we have presented.

\section*{Acknowledgements} 
The author would like to thank Karen Collins for her encouragement, insightful conversations, and suggestions relating to this work. Thanks are also due to an anonymous referee for helpful suggestions.

\bibliographystyle{plain}
\bibliography{references}
\vspace{1cm}

Department of Mathematics, Oberlin College, Oberlin Ohio\

igossett@oberlin.edu\

\end{document}